\documentclass[12pt]{amsart}

\usepackage{amsmath, amsthm, amscd, amsfonts, amssymb, graphicx, color}
\usepackage{times}
\usepackage{bm}
\usepackage{natbib}
\newtheorem{theorem}{Theorem}

\def\ca{\mathring}

\begin{document}

%\jname{Biometrika}
%%% The year, volume, and number are determined on publication
%\jyear{}
%\jvol{}
%\jnum{}
%%% The \doi{...} and \accessdate commands are used by the production team
%%\doi{10.1093/biomet/asm023}
%%\accessdate{Advance Access publication on 31 July 2008}
%\copyrightinfo{\Copyright\ 2011 Biometrika Trust\goodbreak {\em Printed in Great Britain}}

%% These dates are usually set by the production team
%\received{January 2008}
%\revised{July 2008}

%% The left and right page headers are defined here:
%\markboth{THOMAS J. DICICCIO, TODD A. KUFFNER \and G. ALASTAIR YOUNG}{Objective Bayes and the signed root
%statistic}

\title{Objective Bayes, conditional inference and the signed root likelihood ratio statistic}

\author{THOMAS J. DICICCIO}
\address{Department of Social Statistics,
Cornell University, Ithaca, New York 14853, U.S.A.} \email{tjd9@cornell.edu}

\author{TODD A. KUFFNER}
\address{Department of Mathematics, American University of Beirut, Beirut 1107 2020, Lebanon.} \email{tk30@aub.edu.lb}

\author{G. ALASTAIR YOUNG}
\address{Department of Mathematics, Imperial College London, London SW7 2AZ,
U.K.} \email{alastair.young@imperial.ac.uk}

\begin{abstract}
Bayesian properties of the signed root likelihood ratio statistic are analysed. Conditions for first-order probability matching are derived by examination of the Bayesian posterior and frequentist means of this statistic. Second-order matching conditions are shown to arise from matching of the Bayesian posterior and frequentist variances of a mean-adjusted version of the signed root statistic. Conditions for conditional probability matching in ancillary statistic models are derived and discussed.
\end{abstract}


\maketitle

\noindent\section{Introduction}

In problems concerning inference on a scalar interest parameter in the presence of a nuisance parameter, the signed root likelihood ratio statistic is a fundamental object of statistical methodology.  The focus of this paper is an analysis of the signed root statistic from an objective Bayes perspective, where the issue of identification of prior distributions which display the property of probability matching is central. Under probability matching, quantiles of the Bayesian posterior distribution have the property of being confidence limits in the conventional, frequentist sense.

Considering inference based on a sample of size $n$, we establish a number of results. If we choose the prior distribution in a Bayesian analysis so that the frequentist and Bayesian posterior means of the signed root statistic match to $O_p(n^{-1})$, then the prior achieves first-order probability matching: the Bayesian $1-\alpha$ posterior quantile has frequentist coverage $1-\alpha+O(n^{-1})$. It is observed that such matching of frequentist and Bayesian posterior means occurs when the prior satisfies the conditions noted by Peers (1965), who extended to the nuisance parameter case work of Welch \& Peers (1963). We further obtain a simple condition, related to the Bayesian posterior and frequentist variances of a mean-adjusted version of the signed root statistic, under which the prior is second-order probability matching, so that the frequentist coverage of the Bayesian $1-\alpha$ quantile is $1-\alpha+O(n^{-3/2})$. This condition is shown to be equivalent to that established by Mukerjee \& Dey (1993) and Mukerjee \& Ghosh (1997) by an analytically more elaborate shrinkage argument. Our analysis therefore provides a transparent and intuitive interpretation, in terms of the distributional properties of the signed root statistic, for existing conditions for probability matching priors.

In particular statistical problems, specifically those involving inference on canonical parameters in multi-parameter exponential family models and in models admitting ancillary statistics, the appropriate frequentist inference is a conditional one, so that the relevant objective Bayesian notion is that of a conditional probability matching prior. We further provide an analysis of the conditional frequentist behaviour of Bayesian posterior quantiles in ancillary statistic models, as was carried out for the exponential family context by DiCiccio \& Young (2010). Their approach involves matching of higher-order conditional frequentist and Bayesian asymptotics, yielding simple conditions for probability matching. In the ancillary statistic context, this approach does not, however, yield any general, easily interpretable conditions, unlike the methodologies applied here. By considering the distributional properties of the signed root statistic, we note first that if the unconditional first-order probability matching condition of Peers (1965), which we refer to as the Welch--Peers condition, is satisfied, then the prior automatically enjoys the property of first-order conditional probability matching. We then establish our key result, which identifies a simple condition involving the Bayesian posterior and conditional frequentist means of the signed root statistic, under which the conditions for second-order probability matching in an unconditional sense ensure also second-order conditional probability matching.

%Section 2 provides preliminaries and notation for our discussion. Section 3 examines the mean of the signed root %statistic and establishes by this analysis that unconditional first-order probability matching is obtained when %the frequentist and Bayesian means match to $O_p(n^{-1})$, which happens when the Welch-Peers condition holds. %This provides another route, additional to the existing objective Bayesian literature, as considered for example %by Datta \& Mukerjee (2004), to identification of probability matching priors. Section 4 considers
%a higher-order analysis of the Bayesian mean of the signed root statistic, which facilitates our analysis in %Section 5 of second-order probability matching priors, via consideration of the frequentist and Bayesian %posterior variances of a mean-adjusted version of the statistic. Section 6 considers conditional inference. It %provides further discussion of the variance of the mean-adjusted statistic and provides proof of our main result %concerning conditional probability matching priors. Concluding discussion and illustration of the results is %given in Section 7.

\bigskip\noindent
\section{Notation}

Consider a random vector $Y=(Y_1,\ldots,Y_n)$ having continuous probability distribution that depends on an unknown $(q+1)$-dimensional parameter $\theta=(\theta^1,\ldots,\theta^{q+1})$, and denote the log likelihood function for $\theta$ based on $Y$ by $L(\theta)$. Suppose that $\theta$ is partitioned in the form $\theta=(\psi,\phi)$, where $\psi$ is a scalar interest parameter and $\phi$ is a $q$-dimensional nuisance parameter. Let $\hat\theta=(\hat\psi,\hat\phi)$ be the overall maximum likelihood estimator of $\theta$, and let $\tilde\theta(\psi)=\{\psi,\tilde\phi(\psi)\}$ be the constrained maximum likelihood estimator of $\theta$ for given $\psi$.
The log profile likelihood function for $\psi$ is $M(\psi)=L\{\tilde\theta(\psi)\}$.

In the asymptotic calculations that follow, standard conventions for denoting arrays and summation are used. For these conventions, it is understood that indices $i, j, k, \ldots$ range over $2,\ldots,q+1$, and that indices $r,s,t,\ldots$ range over $1,\ldots,q+1$. Summation over the relevant range is implied for any index appearing in an expression both as a subscript and as a superscript.
Differentiation of the functions $L(\theta)$ and $M(\psi)$ is indicated by subscripts, so $L_r(\theta)=\partial L(\theta)/\partial\theta^r$, $L_{rs}(\theta)=\partial^2 L(\theta)/\partial\theta^r \partial\theta^s$, $M_1(\psi)=\partial M(\psi)/\partial\psi$, $M_{11}(\psi)=\partial^2 M(\psi)/\partial\psi^2$, etc.
In this notation, $L_r(\hat\theta)=0$ $(r=1,\ldots,q+1)$ and $ M_1(\hat\psi)=0$.
Evaluation of the derivatives of $L(\theta)$ at $\hat\theta$ and the derivatives of $M(\psi)$ at $\hat\psi$ is indicated by placing a circumflex above the appropriate quantity; for example, $\hat L_{r}=L_{r}(\hat\theta)=0$, $\hat L_{rs}=L_{rs}(\hat\theta)$, $\hat M_1=M_1(\hat\psi)$, $\hat M_{11}=M_{11}(\hat\psi)$, etc.
Let $\lambda_r={ E}\{L_r(\theta)\}=0$, $\lambda_{rs}={ E}\{L_{rs}(\theta)\}$, $\lambda_{rst}={ E}\{L_{rst}(\theta)\}$, etc., and define $l_r=L_r(\theta)-\lambda_r=L_r(\theta)$, $l_{rs}=L_{rs}(\theta)-\lambda_{rs}$, $l_{rst}=L_{rst}(\theta)-\lambda_{rst}$, etc. The constants $\lambda_{rs}$, $\lambda_{rst}$, etc.\ are assumed to be of order $O(n)$; the variables $l_r$, $l_{rs}$, $l_{rst}$, etc.\ have expectation 0, and they are assumed to be of order $O_p(n^{1/2})$. The joint cumulants of $l_r, l_{rs}$, etc.\ are assumed to be of order $O(n)$. These assumptions are usually satisfied in situations involving independent
observations.

In subsequent calculations, it is useful to extend the $\lambda$-notation: let
$\lambda_{r,s}={E}\{L_r(\theta)L_s(\theta)\}$, $\lambda_{rs,t}={E}\{L_{rs}(\theta)L_t(\theta)\}$, $\lambda_{r,s,t}={E}\{L_r(\theta)L_s(\theta)L_t(\theta)\}$,
etc. Identities involving the $\lambda$'s can be derived by repeated differentiation of the identity $\int\exp\{L(\theta)\}dy=1$; in particular,
$
\lambda_{rs}+\lambda_{r,s}=0$,
$\lambda_{rst}+\lambda_{rs,t}+\lambda_{rt,s}+\lambda_{st,r}+\lambda_{r,s,t}=0.
$
Differentiation of the definition $\lambda_{rs}=\int L_{rs}(\theta)\exp\{L(\theta)\}dy$ yields
$\lambda_{rs/t}=\lambda_{rst}+\lambda_{rs,t},$
where $\lambda_{rs/t}=\partial\lambda_{rs}/\partial\theta^t$.

Let $(\lambda^{rs}), (L^{rs})$ and $(\hat L^{rs})$ be the $(q+1)\times(q+1)$ matrix inverses of $(\lambda_{rs}), (L_{rs})$ and $(\hat L_{rs})$, respectively.
Define $\tau^{rs}=\lambda^{r1}\lambda^{s1}/\lambda^{11}$, $\hat T^{rs}=\hat L^{r1}\hat L^{s1}/\hat L^{11}$, $\nu^{rs}=\lambda^{rs}-\tau^{rs}$, and $\hat V^{rs}=\hat L^{rs}-\hat T^{rs}$.
Note that $\lambda^{rs}$, $\tau^{rs}$, and $\nu^{rs}$ are all of order $O(n^{-1})$, and $\hat L^{rs}$, $\hat T^{rs}$, and $\hat V^{rs}$ are all of order $O_p(n^{-1})$.
Furthermore, $1/\lambda^{11}=\lambda_{r1}\lambda_{s1}\tau^{rs}$ is of order $O(n)$ and $1/\hat L^{11}=\hat L_{r1}\hat L_{s1}\hat T^{rs}$ is of order $O_p(n)$.
Note that $\tau^{r1}=\lambda^{r1}$ and $\nu^{r1}=0$; thus, the entries of $(q+1)\times(q+1)$ matrices $(\nu^{rs})$ and $(\hat V^{rs})$ are all 0 except for the lower right-hand submatrices $(\nu^{ij})$ and $(\hat V^{ij})$, which are the inverses of $(\lambda_{ij})$ and $(\hat L_{ij})$, respectively.

\bigskip
\noindent \section{Signed root statistic and probability matching}

The likelihood ratio statistic $W(\psi)=2\{M(\hat\psi)-M(\psi)\}$ is useful for testing the null hypothesis $H_0: \psi=\psi_0$ against the two-sided alternative $H_a:\psi\ne\psi_0$ or for constructing two-sided confidence intervals. However, for a scalar interest parameter, it is natural to conduct one-sided tests, where the alternative is either $H_a:\psi>\psi_0$ or $H_a:\psi<\psi_0$, or to construct one-sided, upper or lower, confidence limits.  One-sided inferences can be achieved by considering the signed square root of the likelihood ratio statistic $R(\psi)={\rm sgn(\hat\psi-\psi)}\{W(\psi)\}^{1/2}$,
%The likelihood ratio statistic $W(\psi)$ typically has the $\chi^2_1$ distribution to error of order $O(n^{-1})$, %so the $\chi^2_1$ approximation produces confidence regions having coverage error of that order. The order of %error in the $\chi^2_1$ approximation to the distribution of $W(\psi)$ can be reduced to $O(n^{-2})$ by Bartlett %correction. In contrast, the signed root $R(\psi)$
which has the standard normal distribution, $N(0,1)$, to error of order $O(n^{-1/2})$, so that the $N(0,1)$ approximation produces one-sided confidence limits having coverage error of that order. The order of error in the $N(0,1)$ approximation to the distribution of $R(\psi)$ can be reduced (DiCiccio \& Stern, 1994a) to $O(n^{-1})$ by correcting for the mean of $R(\psi)$.
% if $\mu(\theta)$ denotes the mean of $R(\psi)$, then $\mu(\theta)$ is of order $O(n^{-1/2})$, and the %mean-adjusted statistic $R(\psi)-\mu(\theta)$ has the $N(0,1)$ distribution to error of order $O(n^{-1})$. This %order of error holds if the mean $\mu(\theta)$ is approximated by an estimate $\breve\mu$ such that %$\breve\mu=\mu(\theta)+O_p(n^{-1})$. Examples include $R(\psi)-\mu\{\tilde\theta(\psi)\}$ and %$R(\psi)-\mu(\hat\theta)$. Applying the $N(0,1)$ approximation to the distribution of $R(\psi)-\breve\mu$ yields %one-sided confidence limits having coverage error of order $O(n^{-1})$.

The signed root $R(\psi)$ behaves (DiCiccio \& Stern, 1994b) identically from a Bayesian perspective: the posterior distribution of $R(\psi)$ is standard normal to error of order $O(n^{-1/2})$, and the order of error in the $N(0,1)$ approximation to the posterior distribution of $R(\psi)$ can be reduced to $O(n^{-1})$ by correcting for the posterior mean. These order statements are asserted conditionally given the data, so it is appropriate to use $O(\cdot)$ rather than $O_p(\cdot)$ to describe the errors associated with the $N(0,1)$ approximation. To distinguish frequentist and Bayesian inference, let $\mu_F=\mu_F(\theta)$ denote the frequentist mean of $R(\psi)$ and let $\mu_B=\mu_B(Y)$ denote the posterior mean. Applying the $N(0,1)$ approximation to the posterior distribution of $R(\psi)-\mu_B$ yields approximate posterior percentage points for $\psi$ having posterior probability error of order $O(n^{-1})$. If the prior distribution is chosen so that $\mu_B=\mu_F+O_p(n^{-1})$, then the posterior distribution of $R(\psi)-\mu_B$ coincides with the frequentist one to error of order $O(n^{-1})$. Thus, for such a prior distribution, the upper $1-\alpha$ posterior percentage point is necessarily an approximate upper $1-\alpha$ frequentist confidence limit having coverage error of order $O(n^{-1})$. As noted, Peers (1965) derived a condition that the prior distribution must satisfy in order for it to have this first-order probability matching property, although not by the method considered here of matching the Bayesian posterior and frequentist means of the signed root statistic.

DiCiccio \& Stern (1994a) showed that $\mu_F$ is
$$
\mu_F=-\textstyle{1\over 2}\eta\lambda_{rst}\lambda^{r1}\lambda^{st}
-\textstyle{1\over 6}\eta^3\lambda_{rst}\lambda^{r1}\lambda^{s1}\lambda^{t1}
+\eta\lambda_{rs/t}\lambda^{r1}\lambda^{st}
+\textstyle{1\over 2}\eta^3\lambda_{rs/t}\lambda^{r1}\lambda^{s1}\lambda^{t1}+O(n^{-3/2}),
$$
where $\eta=(-\lambda^{11})^{-1/2}$. A general expression for $\mu_B$ is derived in the Appendix.
%\begin{equation}
%\mu_B=-{\textstyle{1\over 2}}\hat H\hat L_{rst}\hat L^{r1}\hat L^{rs}
%-{\textstyle{1\over 6}}\hat H^3\hat L_{rst}\hat L^{r1}\hat L^{s1}\hat L^{t1}
%+\hat H{{\partial\log\pi(\theta)}\over{\partial\theta^r}}\bigg|_{\theta=\hat\theta}\hat L^{r1}+O(n^{-3/2}).
%\end{equation}

To compare $\mu_B$ and $\mu_F$, note that
$$\mu_B=-{\textstyle{1\over 2}}\eta\lambda_{rst}\lambda^{r1}\lambda^{rs}
-{\textstyle{1\over 6}}\eta^3\lambda_{rst}\lambda^{r1}\lambda^{s1}\lambda^{t1}
+\eta{{\partial\log\pi(\theta)}\over{\partial\theta^r}}\lambda^{r1}+O_p(n^{-1})$$
in the frequentist sense; hence, $O_p(\cdot)$ is used for the error term in place of $O(\cdot)$. Thus, the condition $\mu_B=\mu_F+O_p(n^{-1})$ is met when the prior satisfies
$$\eta{{\partial\log\pi(\theta)}\over{\partial\theta^r}}\lambda^{r1}=
\eta\lambda_{rs/t}\lambda^{r1}\lambda^{st}
+\textstyle{1\over 2}\eta^3\lambda_{rs/t}\lambda^{r1}\lambda^{s1}\lambda^{t1}.
$$
A standard result of matrix algebra gives $(\lambda^{uv})_{/t}=-\lambda_{rs/t}\lambda^{ru}\lambda^{sv}$, so it follows that
$$
\sum_{r}{{\partial\lambda^{r1}}\over{\partial\theta^r}}=-\lambda_{rs/t}\lambda^{r1}\lambda^{st},\quad
{{\partial\eta}\over{\partial\theta^r}}=-\textstyle{1\over 2}\eta^3\lambda_{st/r}\lambda^{s1}\lambda^{t1};
$$
consequently, the condition on the prior can be written as
\begin{equation}
\eta{{\partial\log\pi(\theta)}\over{\partial\theta^r}}\lambda^{r1}
=-\sum_r{{\partial(\eta\lambda^{r1})}\over{\partial\theta^r}},
\end{equation}
which is what we have termed the Welch--Peers condition.

%When the interest parameter is orthogonal to the nuisance parameters, i.e., when $\lambda^{i1}=0$, then %$\lambda^{11}=1/\lambda_{11}$ and $\eta=(-\lambda_{11})^{1/2}$, so the Welch-Peers condition reduces to
%$$(-\lambda_{11})^{-1/2} {{\partial\log\pi(\theta)}\over{\partial\psi}}
%=-{{\partial (-\lambda_{11})^{-1/2}}\over{\partial\psi}},$$
%which yields $\pi(\theta)\propto\sqrt{-\lambda_{11}}g(\phi),$ where $g(\phi)$ is an arbitrary function of the %nuisance parameter $\phi$: see Tibshirani (1989).

In many cases, Welch--Peers priors are not unique; for example, under parameter orthogonality, the prior is (Tibshirani, 1989) essentially the Jeffreys prior on the interest parameter, multiplied by an arbitrary function of the nuisance parameter. As a means to choose among the Welch--Peers priors, it is natural to attempt to determine those that are second-order probability matching, i.e., those for which the percentage points are one-sided confidence limits having coverage error of order $O(n^{-3/2})$. This problem has been well-studied, and Mukerjee \& Dey (1993) and Mukerjee \& Ghosh (1997) have given a differential equation, auxiliarly to the Welch--Peers condition, that the prior must satisfy for it to be second-order probability matching: see Datta \& Mukerjee (2004, Theorem 2.4.1) for a summary. The objective here is to demonstrate that the auxiliary condition can be developed by considering a mean-adjusted version of the signed root statistic.

DiCiccio \& Stern (1993) showed that the posterior expectation of $\{R(\psi)\}^2$ is $1+a_B+O(n^{-3/2})$, where
\begin{eqnarray*}
&&a_B={\textstyle{1\over 4}}(\hat L^{rs}\hat L^{tu}-\hat V^{rs}\hat V^{tu})\hat L_{rstu}
-{\textstyle{1\over 4}}(\hat L^{ru}\hat L^{st}\hat L^{vw}-\hat V^{ru}\hat V^{st}\hat V^{vw})\hat L_{rst}\hat L_{uvw}\\
&& -{\textstyle{1\over 6}}(\hat L^{ru}\hat L^{sw}\hat L^{tv}-\hat V^{ru}\hat V^{sw}\hat V^{tv})\hat L_{rst}\hat L_{uvw}
+(\hat L^{rs}\hat L^{tu}-\hat V^{rs}\hat V^{tu})\hat L_{rst}\hat\Pi_u
-(\hat L^{rs}-\hat V^{rs})\hat\Pi_{rs},
\end{eqnarray*}
and $\hat\Pi_r=\Pi_r(\hat\theta)$,  $\hat\Pi_{rs}=\Pi_{rs}(\hat\theta)$, with $\Pi_r(\theta)=\pi_r(\theta)/\pi(\theta)$, $\Pi_{rs}(\theta)=\pi_{rs}(\theta)/\pi(\theta)$, $\pi_r(\theta)=\partial\pi(\theta)/\partial\theta^r$, $\pi_{rs}(\theta)=\partial^2\pi(\theta)/\partial\theta^r\partial\theta^s$. It follows that the posterior variance of $R(\psi)-\mu_B$ is
$$\sigma^2_B=1+a_B-\mu_B^2+O(n^{-3/2}).$$ When the prior satisfies condition (1), from a frequentist perspective, since $\mu_B$ is of order $O_p(n^{-1/2})$  and $\mu_B=\mu_F+O_p(n^{-1})$, where $\mu_F$ is of order $O(n^{-1/2})$, it follows that $\mu^2_B=\mu^2_F+O_p(n^{-3/2})$, where $\mu_F^2$ is of order $O(n^{-1})$; hence, the posterior variance satisfies
$$\sigma^2_B=1+a_B-\mu_F^2+O_p(n^{-3/2}).$$

DiCiccio \& Stern (1994b) showed that the frequentist expectation of $\{R(\psi)\}^2$ is $1+a_F+O(n^{-3/2})$, where
\begin{eqnarray*}
a_F&=&(\lambda^{rs}\lambda^{tu}-\nu^{rs}\nu^{tu})({\textstyle{1\over 4}}\lambda_{rstu}-\lambda_{rst/u}+
\lambda_{rt/su})\\
&&\qquad-(\lambda^{ru}\lambda^{st}\lambda^{vw}-\nu^{ru}\nu^{st}\nu^{vw})({\textstyle{1\over 4}}\lambda_{rst}\lambda_{uvw}-\lambda_{rst}\lambda_{uv/w}+\lambda_{rs/t}\lambda_{uv/w})\\
&&\qquad\qquad-(\lambda^{ru}\lambda^{sw}\lambda^{tv}-\nu^{ru}\nu^{sw}\nu^{tv})({\textstyle{1\over 6}}\lambda_{rst}\lambda_{uvw}-\lambda_{rst}\lambda_{uv/w}+\lambda_{rs/t}\lambda_{uv/w})+O(n^{-3/2}).
\end{eqnarray*}
It is shown in the appendix that ${E}(\mu_B)=\mu_F+O(n^{-3/2})$, and so it follows as a particular case of the results of DiCiccio \& Stern (1994a) that the frequentist variance of $R(\psi)-\mu_B$ is
$$\sigma^2_F=1+a_F+2\eta\mu_{F/r}\lambda^{r1}-\mu_F^2+O(n^{-3/2}),$$
where $\mu_{F/r}=\partial\mu_F/\partial\theta^r$.
%Almost all of the corrected versions of $R(\psi)$  appearing in the literature that improve the accuracy of the %standard normal approximation to the distribution of $R(\psi)$ are of the form $R(\psi)-D(\psi)$ where $D(\psi)$ %has the following properties: $D(\psi)$ and its derivatives taken with respect to $\psi$ are of order %$O_p(n^{-1/2})$; $D(\psi)=\mu_F+O_p(n^{-1})$; and ${ E}\{D(\psi)\}=\mu_F+O(n^{-3/2})$. The correction term %$D(\psi)=\mu_B$ considered here has been shown above to have these properties. In the general context, DiCiccio %\& Stern (1994a) showed that the variance of $R(\psi)-D(\psi)$ is
%$${\rm var}\{R(\psi)-D(\psi)\}=1+a_F+2\eta\mu_{F/r}\lambda^{r1}-\mu_F^2+2\eta^{-1}E\{\partial %D(\psi)/\partial\psi\}+O(n^{-3/2}).$$
%In the particular case $D(\psi)=\mu_B$, we have $\partial D(\psi)/\partial\psi=0$, so $E\{\partial %D(\psi)/\partial\psi\}=0$.

From both Bayesian and frequentist perspectives, the third- and higher-order cumulants of $R(\psi)-\mu_B$ are (DiCiccio \& Stern, 1994a) of order $O(n^{-3/2})$ or smaller. Thus, the marginal distribution of $\{R(\psi)-\mu_B\}/\sigma_F$ and the posterior distribution of $\{R(\psi)-\mu_B\}/\sigma_B$ are both standard normal to error of order $O(n^{-3/2})$; moreover, if $\pi(\theta)$ is a prior density such that $\sigma^2_B=\sigma^2_F+O_p(n^{-3/2})$, then $\pi(\theta)$ is a second-order probability matching prior: a formal proof follows the same argument as that used in our main result in Section 4 below. The condition under which $\pi(\theta)$ will be a second-order probability matching prior is therefore that
$$a_B=a_F+2\eta\mu_{F/r}\lambda^{r1}+O_p(n^{-3/2}).$$
Since
\begin{eqnarray*}
a_B &=&{\textstyle{1\over 4}}(\lambda^{rs}\lambda^{tu}-\nu^{rs}\nu^{tu})\lambda_{rstu}
-{\textstyle{1\over 4}}(\lambda^{ru}\lambda^{st}\lambda^{vw}-\nu^{ru}\nu^{st}\nu^{vw})\lambda_{rst}\lambda_{uvw}\\
&& \qquad-{\textstyle{1\over 6}}(\lambda^{ru}\lambda^{sw}\lambda^{tv}-\nu^{ru}\nu^{sw}\nu^{tv})\lambda_{rst}\lambda_{uvw}
+(\lambda^{rs}\lambda^{tu}-\nu^{rs}\nu^{tu})\lambda_{rst}\Pi_u \\
&&\qquad \qquad-(\lambda^{rs}-\nu^{rs})\Pi_{rs}+O_p(n^{-3/2}),
\end{eqnarray*}
this condition can be expressed as
\begin{eqnarray*}
\lefteqn{
a_F+2\eta\mu_{F/r}\lambda^{r1}={\textstyle{1\over 4}}(\lambda^{rs}\lambda^{tu}-\nu^{rs}\nu^{tu})\lambda_{rstu}
-{\textstyle{1\over 4}}(\lambda^{ru}\lambda^{st}\lambda^{vw}-\nu^{ru}\nu^{st}\nu^{vw})\lambda_{rst}\lambda_{uvw}}\\
&&-{\textstyle{1\over 6}}(\lambda^{ru}\lambda^{sw}\lambda^{tv}-\nu^{ru}\nu^{sw}\nu^{tv})\lambda_{rst}\lambda_{uvw}
+(\lambda^{rs}\lambda^{tu}-\nu^{rs}\nu^{tu})\lambda_{rst}\Pi_u-(\lambda^{rs}-\nu^{rs})\Pi_{rs}.
\end{eqnarray*}
%that is,
%\begin{eqnarray*}
%\lefteqn{-(\lambda^{rs}-\nu^{rs})\Pi_{rs}+(\lambda^{rs}\lambda^{tu}-\nu^{rs}\nu^{tu})
%\lambda_{rst}\Pi_u=} \\
%& &-(\lambda^{rs}\lambda^{tu}-\nu^{rs}\nu^{tu})\{\lambda_{rst/u}-\lambda_{rt/su}\}
%+(\lambda^{ru}\lambda^{st}\lambda^{vw}-\nu^{ru}\nu^{st}\nu^{vw})\{\lambda_{rst}
%\lambda_{uv/w}-\lambda_{rs/t}\lambda_{uv/w}\}\\
%&&+(\lambda^{ru}\lambda^{sw}\lambda^{tv}-\nu^{ru}\nu^{sw}\nu^{tv})\{\lambda_{rst}\lambda_{uv/w}-\lambda_{rs/t}\lambda_{uv/w}\}+2\eta\mu_{F/r}\lambda^{r1}.
%\end{eqnarray*}
By assuming that the prior density $\pi(\theta)$ satisfies the first-order probability-matching condition (1), the condition for second-order probability matching reduces to
%\begin{eqnarray*}
%\lefteqn{
%\tau^{rs}\Pi_{rs}-\tau^{rs}\lambda^{tu}\lambda_{rst}\Pi_u=(\tau^{rs}\lambda^{tu}-
%{\textstyle{2\over %3}}\tau^{rs}\tau^{tu})\lambda_{rst/u}-\lambda^{ru}\tau^{st}\lambda^{vw}\lambda_{rst}\lambda_{uv/w}} \\
%&&-(2\lambda^{ru}\lambda^{sw}\tau^{tv}-2\lambda^{ru}\tau^{sw}\tau^{tv}-\tau^{ru}
%\lambda^{sw}\tau^{tv}+\tau^{ru}\tau^{sw}\tau^{tv})\lambda_{rst}\lambda_{uv/w}\\
%&&+(\tau^{ru}\lambda^{st}\lambda^{vw}-\tau^{ru}\lambda^{st}\tau^{vw}+{\textstyle{1\over %2}}\tau^{ru}\tau^{st}\tau^{vw})\lambda_{rs/t}\lambda_{uv/w} \\
%&&+(\tau^{ru}\lambda^{sw}\lambda^{tv}-\tau^{ru}\lambda^{sw}\tau^{tv})\lambda_{rs/t}\lambda_{uv/w},
%\end{eqnarray*}
%which can be written somewhat more succinctly as
\begin{eqnarray*}
\lefteqn{
\tau^{rs}\Pi_{rs}-\tau^{rs}\lambda^{tu}\lambda_{rst}\Pi_u=(\tau^{rs}\nu^{tu}+
{\textstyle{1\over 3}}\tau^{rs}\tau^{tu})\lambda_{rst/u}-\lambda^{ru}\tau^{st}\lambda^{vw}\lambda_{rst}\lambda_{uv/w}} \\
&&-(\lambda^{ru}\nu^{sw}\tau^{tv}+\nu^{ru}\nu^{sw}\tau^{tv})\lambda_{rst}\lambda_{uv/w}
+(\tau^{ru}\lambda^{st}\nu^{vw}+{\textstyle{1\over 2}}\tau^{ru}\tau^{st}\tau^{vw})\lambda_{rs/t}\lambda_{uv/w}\\
&&+\tau^{ru}\lambda^{sw}\nu^{tv}\lambda_{rs/t}\lambda_{uv/w},
\end{eqnarray*}
which may be written more succinctly as
\begin{equation}
\sum_r{{\partial\{\pi(\nu^{rs}+{\textstyle{1\over 3}}\tau^{rs})\tau^{tu}\lambda_{stu}\}}
\over{\partial\theta^r}}+\sum_{r,s}{{\partial^2(\pi\tau^{rs})}\over{\partial\theta^r\partial\theta^s}}=0.
\end{equation}

This condition is that given by Mukerjee \& Dey (1993), who considered a scalar nuisance parameter, and by Mukerjee \& Ghosh (1997), who considered a vector nuisance parameter.

\smallskip
\noindent\section{Conditional inference}

The main motivation for our analysis is to investigate the conditional frequentist properties of approximate confidence limits obtained from Welch--Peers priors. Suppose that $A$ is an ancillary statistic such that $(\hat\theta,A)$ is sufficient.
To undertake calculations with respect to the conditional distribution of $Y$ given $A$, or equivalently, with respect to the conditional distribution of $(\hat\theta,A)$ given $A$, it is useful to consider the versions of the $\lambda$'s obtained by applying their definitions to the conditional distribution.
The resulting conditional quantities are distinguished notationally from the unconditional ones by the inclusion of the accent symbol ``$\ca{~~~}$''.
Since the conditional log likelihood function $\ca L(\theta)$ differs from the unconditional log likelihood function $L(\theta)$ by a constant, i.e., a quantity that depends on $A$ but not on $\theta$, it follows that ${\ca L}_r=L_r$, ${\ca L}_{rs}=L_{rs}$, etc.
Thus, $\ca\lambda_r={ E}_{Y|A}(L_r)$, $\ca\lambda_{rs}={ E}_{Y|A}(L_{rs})$, etc.
Note that $\ca\lambda_r=0$.
The quantities $\ca\lambda_{rs}$, $\ca\lambda_{rst}$, etc.\ are random variables depending on $A$, and they are assumed to be of order $O_p(n)$.
The variables $\ca l_r=L_r$, $\ca l_{rs}=L_{rs}-\ca\lambda_{rs}$, etc.\ all have zero conditional expectation, and hence they have zero unconditional expectation, and they are assumed to be of order $O_p(n^{1/2})$.
The joint conditional cumulants of $\ca l_r, \ca l_{rs}$, etc.\ depend on $A$, and they are assumed to be of order $O(n)$ given $A$ and of order $O_p(n)$ unconditionally.
The identities that hold for the $\lambda$'s immediately carry over to the $\ca\lambda$'s.

In the calculations that follow, it is necessary to take into account the differences between the $\lambda$'s and the $\ca\lambda$'s. To describe the difference between $\lambda_{rs}$ and $\ca\lambda_{rs}$, first note that
${E}_{A}(\ca\lambda_{rs})={ E}_{A}\{{ E}_{Y|A}(L_{rs})\}={ E}_Y(L_{rs})=\lambda_{rs};$
moreover,
\begin{eqnarray*}
{\rm var}_A(\ca\lambda_{rs})&=&{\rm var}_A\{{ E}_{Y|A}(L_{rs})\}
={\rm var}_Y(L_{rs})-{ E}_A\{{\rm var}_{Y|A}(L_{rs})\}\\
&=&{\rm var}_Y(l_{rs})-{ E}_A\{{\rm var}_{Y|A}(\ca l_{rs})\}
=O(n)-{\ E}_A\{O_p(n)\}
=O(n).
\end{eqnarray*}
With respect to the distribution of $A$, $\ca\lambda_{rs}$ has mean $\lambda_{rs}$ and variance of order $O(n)$; thus, $\ca\lambda_{rs}-\lambda_{rs}$,  $\ca\lambda_{rst}-\lambda_{rst}$, etc.\ have expectation 0 and variance of order $O(n)$. It is assumed that these quantities are of order $O_p(n^{1/2})$ and have joint cumulants of order $O(n)$ with respect to the distribution of $A$. In particular, $\ca\lambda_{rs}=\lambda_{rs}+O_p(n^{1/2})$, $\ca\lambda_{rst}=\lambda_{rst}+O_p(n^{1/2})$, etc.

Assume that differentiation of the identity $\ca\lambda_{rs}=\lambda_{rs}+O_p(n^{1/2})$ yields $\ca\lambda_{rs/t}=\lambda_{rs/t}+O_p(n^{1/2})$.  We note that, as a rule, differentiation of an asymptotic
relation will preserve the asymptotic order, but that care is necessary; see Barndorff-Nielsen \& Cox (1994, Exercise 5.4).  The difference between $\ca\lambda_{rs/t}$ and $\lambda_{rs/t}$ indicated here actually constitutes an additional assumption of our calculations, though validity is immediate in particular cases, such as the example considered in Section 5. Then,
$$
\ca\lambda_{rs,t}=\ca\lambda_{rs/t}-\ca\lambda_{rst}
=\lambda_{rs/t}-\lambda_{rst}+O_p(n^{1/2})
=\lambda_{rs,t}+O_p(n^{1/2}).
$$

By working with the conditional density of $Y$ given $A$ in place of the marginal density of $Y$, it follows that $R(\psi)$ is conditionally $N(0,1)$ to error of order $O(n^{-1/2})$ and that the error in the $N(0,1)$ approximation to the conditional distribution of $R(\psi)$ can be reduced to order $O(n^{-1})$ by adjusting for the conditional mean of $R(\psi)$. Denote the conditional mean by $\ca\mu_F$; then $R(\psi)-\ca\mu_F$ has conditionally the $N(0,1)$ distribution  to error of order $O(n^{-1})$.

The calculations of DiCiccio \& Stern (1994a) can be applied to the conditional distribution of $R(\psi)$ to show that
$$
\ca\mu_F=-\textstyle{1\over 2}\ca\eta\ca\lambda_{rst}\ca\lambda^{r1}\ca\lambda^{st}
-\textstyle{1\over 6}\ca\eta^3\ca\lambda_{rst}\ca\lambda^{r1}\ca\lambda^{s1}\ca\lambda^{t1}
+\ca\eta\ca\lambda_{rs/t}\ca\lambda^{r1}\ca\lambda^{st}
+\textstyle{1\over 2}\ca\eta^3\ca\lambda_{rs/t}\ca\lambda^{r1}\ca\lambda^{s1}\ca\lambda^{t1}+O(n^{-3/2}),
$$
where $\ca\eta=(-\ca\lambda^{11})^{-1/2}$. Furthermore, the preceding comparisons of the $\lambda$'s and the corresponding $\ca\lambda$'s shows that $\ca\mu_F=\mu_F+O_p(n^{-1})$, and hence, $\mu_B=\ca\mu_F+O_p(n^{-1})$, provided the condition (1) holds. It follows that Welch--Peers priors satisfying (1) produce approximate confidence limits having conditional coverage error of order $O(n^{-1})$. This result is implicit in DiCiccio \& Martin (1993), who compare the Bayesian percentage points under an objective prior with the approximate, conditional confidence limits derived from Barndorff-Nielsen's $r^*$ statistic (Barndorff-Nielsen, 1986); see also Nicolau (1993).

Since $\mu_B=\ca\mu_F+O_p(n^{-1})$, it follows that ${ E}_{Y|A}(\mu_B)=\ca\mu_F+O(n^{-1})$. In resolving conditional properties of second-order probability matching priors, the crucial issue turns out to be to examine circumstances under which ${E}_{Y|A}(\mu_B)=\ca\mu_F+O(n^{-3/2})$ holds. Since $\ca\lambda_r=0$, an argument analogous to that given in the Appendix for $f(\theta)$ defined in (A.1) shows that ${E}_{Y|A}(\mu_B)={E}_{Y|A}\{f(\theta)\}+O(n^{-3/2})$, so the crucial criterion reduces to ${E}_{Y|A}\{f(\theta)\}=\ca\mu_F+O(n^{-3/2})$,
% In the same way as was argued previously,
%\begin{eqnarray*}
%{ E}_{Y|A}\biggl\{H{{\partial\log\pi(\theta)}\over{\partial\theta^r}} L^{r1}\biggr\}&=&
%\ca\eta{{\partial\log\pi(\theta)}\over{\partial\theta^r}}\ca\lambda^{r1}+O(n^{-3/2})\\
%{ E}_{Y|A}(-{\textstyle{1\over 2}} H L_{rst} L^{r1} L^{rs})&=&
%-{\textstyle{1\over 2}} \ca\eta \ca\lambda_{rst} \ca\lambda^{r1} \ca\lambda^{rs}+O(n^{-3/2}),\\
%{ E}_{Y|A}(-{\textstyle{1\over 6}} H^3 L_{rst} L^{r1} L^{s1} L^{t1})&=&
%-{\textstyle{1\over 6}} \ca\eta^3 \ca\lambda_{rst} \ca\lambda^{r1} \ca\lambda^{s1} \ca\lambda^{t1}+O(n^{-3/2}).
%\end{eqnarray*}
%Consequently, ${ E}_{Y|A}\{f(\theta)\}=\ca\mu_F+O(n^{-3/2})$ holds provided
which holds provided
\begin{equation}
\ca\eta{{\partial\log\pi(\theta)}\over{\partial\theta^r}}\ca\lambda^{r1}=\ca\eta\ca\lambda_{rs/t}\ca\lambda^{r1}\ca\lambda^{st}
+\textstyle{1\over 2}\ca\eta^3\ca\lambda_{rs/t}\ca\lambda^{r1}\ca\lambda^{s1}\ca\lambda^{t1}+O(n^{-3/2}).
\end{equation}

Suppose that $\pi$ is a first-order probability-matching prior, so that (1) holds, that satisfies further the condition ${E}_{Y|A}(\mu_B)=\ca\mu_F+O(n^{-3/2})$. Let $\ca\sigma^2_F$ denote the conditional frequentist variance of $R(\psi)-\mu_B$. From arguments similar to the ones given previously that showed $\ca\mu_F=\mu_F+O_p(n^{-1})$, it follows that $\ca\sigma^2_F=\sigma^2_F+O_p(n^{-3/2})$. To be specific, recall that
$$\sigma^2_F=1+a_F+2\eta\mu_{F/r}\lambda^{r1}-\mu_F^2+O(n^{-3/2}),$$
where $a_F+2\eta\mu_{F/r}\lambda^{r1}-\mu_F^2$ is of order $O(n^{-1})$ and can be expressed, to error of order $O(n^{-3/2})$, as a function of the $\lambda$'s. By applying identical calculations, which we note require the condition ${E}_{Y|A}(\mu_B)=\ca\mu_F+O(n^{-3/2})$, to the conditional distribution, it follows that
$$\ca\sigma^2_F=1+\ca a_F+2\ca \eta\ca\mu_{F/r}\ca\lambda^{r1}-\ca\mu_F^2+O(n^{-3/2}),$$ where
$\ca a_F+2\ca \eta\ca\mu_{F/r}\ca\lambda^{r1}-\ca\mu_F^2$ is of order $O_p(n^{-1})$ and can be expressed, to error of order $O_p(n^{-3/2})$, as the identical function as can its unconditional version, with each $\lambda$ being replaced by its corresponding $\ca\lambda$. Since by assumption each $\ca\lambda$ differs from its corresponding $\lambda$ by $O_p(n^{-1/2})$, it follows that $\ca\sigma^2_F=\sigma^2_F+O_p(n^{-3/2})$.
Hence, the condition that ensures $\pi(\theta)$ is a second-order probability-matching prior in the marginal frequentist sense also ensures that it is a second-order probability-matching prior in the conditional frequentist sense.

Thus, if $\pi(\theta)$ is a second-order probability-matching prior in the marginal frequentist sense, then it is also a second-order probability-matching prior in the conditional frequentist sense provided ${E}_{Y|A}(\mu_B)=\ca\mu_F+O(n^{-3/2})$, i.e., provided (3) holds.
%$$
%\ca\eta{{\partial\log\pi(\theta)}\over{\partial\theta^r}}\ca\lambda^{r1}=\ca\eta\ca\lambda_{rs/t}\ca\lambda^{r1}\ca\lambda^{st}
%+\textstyle{1\over 2}\ca\eta^3\ca\lambda_{rs/t}\ca\lambda^{r1}\ca\lambda^{s1}\ca\lambda^{t1}+O(n^{-3/2}).
%$$
This is satisfied if
\begin{equation}
{{\partial\log\pi(\theta)}\over{\partial\theta^r}}\ca\lambda^{r1}=\ca\lambda_{rs/t}\ca\lambda^{r1}\ca\lambda^{st}
+\textstyle{1\over 2}\ca\eta^2\ca\lambda_{rs/t}\ca\lambda^{r1}\ca\lambda^{s1}\ca\lambda^{t1}.
\end{equation}

We summarize our conclusions in the following theorem.

\begin{theorem} Suppose that the prior $\pi$ is such that (4) holds, so that $E_{Y|A}(\mu_{B})=\mathring{\mu}_{F}+O(n^{-3/2})$.  If the prior also satisfies (1) and (2), so that $\sigma_{B}^{2}=\mathring{\sigma}_{F}^{2}+O(n^{-3/2})$, then the Bayesian quantile is second-order conditional probability matching.
\end{theorem}
\emph{Proof}:
Given $Y$, let $\psi_{l} \equiv \psi^{1-\alpha}(\pi,Y)$ be the posterior $1-\alpha$ quantile for the interest parameter, so that ${\rm pr}\{\psi \leq \psi^{1-\alpha}(\pi,Y)|Y\}=1-\alpha$.  Also, given $Y$ and provided the log likelihood is unimodal, $R(\psi)$ is a monotonic decreasing function of $\psi$.  Therefore,
${\rm pr}\{R(\psi) \geq R(\psi_{l})|Y\}=1-\alpha,
$
where $R(\psi)$ is the signed root statistic constructed from $Y$.  That is,

$$
{\rm pr}\Big\{\frac{R(\psi)-\mu_{B}}{\sigma_{B}} \geq \frac{R(\psi_{l})-\mu_B}{\sigma_{B}}|Y\Big\}=1-\alpha.
$$
Since the posterior distribution of $\{R(\psi)-\mu_{B}\}/\sigma_{B}$ is $N(0,1)$ to error of order $O_p(n^{-3/2})$, then by the delta method from Section 2.7 of Hall (1992),
$
{\rm pr}[N(0,1) \geq \{R(\psi_{l})-\mu_{B}\}/{\sigma_{B}}]=1-\alpha+O(n^{-3/2}),
$
so that
$
 \{R(\psi_{l})-\mu_{B}\}/{\sigma_{B}}=z_\alpha +O(n^{-3/2}),
 $
 in terms of the $N(0,1)$ quantile $z_\alpha$ defined by $\Phi(z_\alpha)=\alpha$.

By the monotonicity of $R(\psi)$ given $Y$, the event $\psi \leq \psi_{l}$ is equivalent to the event $R(\psi) \geq R(\psi_{l})$. Thus, from a conditional frequentist perspective, given an ancillary statistic $A=a$, we have, again using the delta method, and by the frequentist distributional result for the mean-adjusted signed root statistic,

\begin{eqnarray*}
{\rm pr}(\psi \leq \psi_{l}|A=a)& =& {\rm pr}\Big\{\frac{R(\psi)-\mu_{B}}{\sigma_{B}}\geq \frac{R(\psi_{l})-\mu_{B}}{\sigma_{B}}|A=a\Big\} \\
& = &{\rm pr}\Big\{\frac{R(\psi)-\mu_{B}}{\mathring{\sigma}_{F}}+O_{p}(n^{-3/2}) \geq \frac{R(\psi_{l})-\mu_{B}}{\sigma_{B}}|A=a\Big\} \\
&=& {\rm pr}\Big\{\frac{R(\psi)-\mu_{B}}{\mathring{\sigma}_{F}} \geq z_\alpha +O_p(n^{-3/2})|A=a\Big\}+O(n^{-3/2}) \\
&=& {\rm pr}\Big\{N(0,1) \geq z_\alpha \Big\}+O(n^{-3/2})
= 1-\alpha+O(n^{-3/2}).
\end{eqnarray*}

%The first equality follows since $R(\psi)$ is monotonically decreasing in $\psi$, the second equality simply %adjusts by the variance, the third equality follows since $\sigma_{B}=\mathring{\sigma}_{F}+O(n^{-3/2})$, the %fourth equality follows by the delta method, the fifth equality follows by the frequentist distributional result %for the mean-adjusted
%signed root statistic combined with the delta method and the last equality is just a re-expression of the %preceding line.

\bigskip
\section{Discussion}

Conditions under which a Bayesian prior on a scalar interest parameter in the presence of a nuisance parameter achieves probability matching have been shown in this paper to have direct interpretation in terms of the frequentist and Bayesian distributional properties of the signed root likelihood ratio statistic. A prior which is first-order probability matching in a marginal sense is necessarily first-order conditional probability matching. A prior which is second-order probability matching in the marginal sense yields second-order conditional probability matching provided a simple condition (4), which may be interpreted as a conditional version of the marginal Welch--Peers condition (1), holds. Second-order unconditional matching is seen
to correspond to matching of the frequentist and Bayesian variances of a specific mean-adjusted version of the signed
root statistic, where the adjustment is by the Bayesian mean, under a prior which ensures first-order matching. Similarly, conditional probability matching is typically
achieved only under a very particular prior specification. These conclusions indicate, to our mind, that second-order probability matching is too stringent a criterion to be useful in practice.
%We argue that more flexible and useful versions of objective Bayes priors might be derived by relaxing the %demands of second-order matching, but discriminating between
%the statistical properties, not related solely to frequentist coverage matching, of competing first-order
%priors.

For illustration, consider a location-scale model, with $Y_1, \ldots, Y_n$ an independent sample from a density of the form $f(y;\mu, \sigma)={\sigma}^{-1}g\{\sigma^{-1}(y-\mu)\}.$
The appropriate conditioning ancillary is the configuration statistic $
A=(A_1, \ldots, A_n)=\{(Y_1-\hat \mu)/{\hat \sigma}, \ldots, (Y_n-\hat \mu)/{\hat \sigma}\},$
and the log-likelihood is of the form
\[
L(\mu, \sigma)=-n\log \sigma+\sum_{i=1}^n h\left(\frac{\hat \mu-\mu}{\sigma}+A_i \frac{\hat \sigma}{\sigma}\right),
\]
with $h(\cdot)=\log g(\cdot)$. Then, since the conditional distribution of $\{(\hat \mu-\mu)/{\sigma}, \hat \sigma /\sigma\}$, given $A=a=(a_1, \ldots, a_n)$, does not depend on $(\mu, \sigma)$, simple calculations show
$\ca \lambda_{\mu \mu}=B/\sigma^2, \ca \lambda_{\mu \sigma}=C/\sigma^2, \ca \lambda_{\sigma \sigma}=D/\sigma^2$, where $B, C, D$ are non-zero functions of $a$. The corresponding unconditional quantities are of the form $ \lambda_{\mu \mu}=B_g/\sigma^2, \lambda_{\mu \sigma}=C_g/\sigma^2, \lambda_{\sigma \sigma}=D_g/\sigma^2$, for constants $B_g, C_g$ and $D_g$ depending only on $g$ and $n$. In the notation of Section 4, it is then immediate that differentiation of the relation $\ca\lambda_{rs}=\lambda_{rs}+O_p(n^{1/2})$ yields the assumed relation $\ca\lambda_{rs/t}=\lambda_{rs/t}+O_p(n^{1/2})$ in this case.

Consider the case where the location parameter $\mu$ is the interest parameter, with $\sigma$ as nuisance parameter. Then the right hand side of the matching condition (4) is easily seen after some algebra to reduce to $\sigma C/E$, where $E=BD-C^2$. The left hand side of (4) is, after some manipulation,
\[
\frac{\sigma^2 D}{E}\frac{\partial \log \pi}{\partial \mu} -\frac{\sigma^2 C}{E} \frac{\partial \log \pi}{\partial \sigma},
\]
so that the matching condition is satisfied if and only if
\[
\frac{\partial \log \pi}{\partial \mu}=0, \quad \frac{\partial \log \pi}{\partial \sigma}=-\frac{1}{\sigma},
\]
which gives $\pi(\mu,\sigma) \propto 1/\sigma$. This prior is easily seen to satisfy (1) and (2), and is therefore second-order conditional probability matching. In fact, the prior is (Lawless, 1982, Appendix G) exact conditional probability matching, and is (Datta \& Mukerjee, 2004, Section 2.5.2) the unique second-order unconditional matching prior.

Analogous calculations yield the same conclusion in the case where the scale parameter $\sigma$ is the interest parameter, with $\mu$ nuisance: the unique solution to the conditional matching condition (4) is $\pi(\mu,\sigma) \propto 1/\sigma$, which is again exact conditional probability matching. Now, however, second-order marginal matching priors are not necessarily unique. In the Cauchy location-scale model, for example, any prior of the form $\pi(\mu,\sigma) \propto d(\mu)/\sigma$, for any smooth positive function $d(\cdot)$, is second-order unconditional matching: see Datta \& Mukerjee (2004, Section 2.5.2).

In general, therefore, second-order conditional probability matching is only achieved by the exact conditional probability matching prior $\pi(\mu,\sigma) \propto 1/\sigma$. We have noted, however, that first-order conditional probability matching is achieved by any first-order unconditional probability matching prior in the class satisfying the Welch--Peers condition (1).
%It may be of some interest to compare, for finite $n$, the conditional probability matching properties of members %of that class. This can be investigated by simulation from the exact conditional density of $\hat \theta \mid %A=a$ and construction of the Bayesian posterior for each such simulated sample.

\setcounter{equation}{0}
\renewcommand{\theequation}{A\arabic{equation}}

\section*{Appendix}
\subsection*{Derivation of $\mu_B$}

Repeated differentiation of the definition $M(\psi)=L\{\tilde\theta(\psi)\}$ yields
$
M_1(\hat\psi)=0,
M_{11}(\hat\psi)= 1/\hat L^{11},
M_{111}(\hat\psi)=\hat L_{rst}\hat L^{r1}\hat L^{s1}\hat L^{t1}/(\hat L^{11})^3$.
%The derivation of these derivatives requires the identities $\tilde\theta(\hat\psi)=\hat\theta$ and %$\tilde\theta^1_{/1}(\psi)=1$. It also requires repeated differentiation of the identity %$L_i\{\tilde\theta(\psi)\}=0$, which yields %$\tilde\theta^r_{/1}(\psi)=L^{r1}\{\tilde\theta(\psi)\}/L^{11}\{\tilde\theta(\psi)\}$, whence %$L_{ir}\{\tilde\theta(\psi)\}\tilde\theta^r_{/1}(\psi)=0$. In particular, $\tilde\theta^r_{/1}(\hat\psi)=\hat %L^{r1}/\hat L^{11}$ and $\hat L_{ir}\tilde\theta^r_{/1}(\hat\psi)=0$.
Taylor expansion about $\hat\psi$ yields
$$W(\psi)=-\hat M_{11}(\hat\psi-\psi)^2
+{\textstyle{1\over 3}}\hat M_{111}(\hat\psi-\psi)^3+O_p(n^{-1}).$$
Consequently,
$$W(\psi)=\{Z(\psi)\}^2-
{\textstyle{1\over 3}}\hat H^3\hat L_{rst}\hat L^{r1}\hat L^{s1}\hat L^{t1}\{Z(\psi)\}^3+O_p(n^{-1}),$$
where $Z(\psi)=(-\hat M_{11})^{1/2}(\hat\psi-\psi)=\hat H(\hat\psi-\psi)$ and $\hat H=(-\hat M_{11})^{1/2}$.
% Thus, $Z(\psi)$ is the familiar studentized statistic, where the standardization is by $\hat H$, the observed %information analogue of $\eta=(-\lambda^{11})^{-1/2}$.
Since $\hat M_{111}=-\hat H^6\hat L_{rst}\hat L^{r1}\hat L^{s1}\hat L^{t1}$,
the signed root statistic $R(\psi)$ has the expansion
$$
R(\psi)=Z(\psi)
-{\textstyle{1\over 6}}\hat H^3\hat L_{rst}\hat L^{r1}\hat L^{s1}\hat L^{t1}\{Z(\psi)\}^2+O_p(n^{-1}).
$$

The Laplace approximation to the marginal posterior density function of $\psi$ given by Tierney \& Kadane (1986) can be written as
$$\pi_{\psi|Y}(\psi)=c\exp\{B(\psi)+M(\psi)-M(\hat \psi)\}\{1+O(n^{-3/2})\},$$
for values of the argument $\psi$ such that $\psi=\hat\psi+O(n^{-1/2})$,
where $c$ is a normalizing constant, $\pi(\theta)=\pi(\psi,\phi)$ is the prior density, and
$$B(\psi)=-{\textstyle{1\over 2}}\log\bigg\{{{|-L_{ij}(\psi,\tilde\phi_\psi)|}
\over{|-L_{ij}(\hat\psi,\hat\phi)|}}\bigg\}
+\log\bigg\{{{\pi(\psi,\tilde\phi_\psi)}
\over{\pi(\hat\psi,\hat\phi)}}\bigg\}.
$$
Differentiation of $B(\psi)$ yields
$$\hat B_1={\textstyle{1\over 2}}\hat H^2\hat L_{rst}\hat L^{r1}\hat L^{rs}
+{\textstyle{1\over 2}}\hat H^4\hat L_{rst}\hat L^{r1}\hat L^{s1}\hat L^{t1}
-\hat H^2{{\partial\log\pi(\theta)}\over{\partial\theta^r}}\bigg|_{\theta=\hat\theta}\hat L^{r1},$$
which is of order $O(1)$.
By Taylor expansion about $\hat\psi$,
$$\pi_{\psi|Y}(\psi)=(2\pi)^{-1/2}\hat H\exp\{-{\textstyle{1\over 2}}\hat H^2(\hat\psi-\psi)^2\}
\{1-\hat B_1(\hat\psi-\psi)-{\textstyle{1\over 6}}\hat M_{111}(\hat\psi-\psi)^3+O(n^{-1})\},
$$
so the marginal posterior density of $Z(\psi)$ has the expansion
$$\pi_{Z(\psi)|Y}(z)=(2\pi)^{-1/2}e^{-z^2/2}
\{1-\hat H^{-1}\hat B_1z-{\textstyle{1\over 6}}\hat H^{-3}\hat M_{111}z^3+O(n^{-1})\},
$$
from which it follows that
\begin{eqnarray*}
\mu_B&=&-\hat H^{-1}\hat B_1-{\textstyle{1\over 2}}\hat H^{-3}\hat M_{111}
-{\textstyle{1\over 6}}\hat H^3\hat L_{rst}\hat L^{r1}\hat L^{s1}\hat L^{t1}+O(n^{-1})\\
&=&-{\textstyle{1\over 2}}\hat H\hat L_{rst}\hat L^{r1}\hat L^{rs}
-{\textstyle{1\over 6}}\hat H^3\hat L_{rst}\hat L^{r1}\hat L^{s1}\hat L^{t1}
+\hat H{{\partial\log\pi(\theta)}\over{\partial\theta^r}}\bigg|_{\theta=\hat\theta}\hat L^{r1}+O(n^{-1}).
\end{eqnarray*}
A more careful analysis that takes higher-order terms into account, i.e., that includes terms of order $O(n^{-1})$, shows that the error term in the preceding formula is actually $O(n^{-3/2})$.

%To be specific, a higher-order expansion of $R(\psi)$ would involve a term in $\{Z(\psi)\}^3$ with coefficient %that is of order $O(n^{-1})$. Moreover, a higher-order expansion  of the marginal  posterior density of $Z(\psi)$ %would involve terms in $z^4$ and $z^6$, each having coefficients of order $O(n^{-1})$. Now, when calculating the %contribution to the expectation from the $\{Z(\psi)\}$ term in $R(\psi)$, the $z^4$ and $z^6$ terms in the %marginal posterior density of $Z(\psi)$ would yield 0. When calculating the contribution to the expectation from %the $\{Z(\psi)\}^2$ term in $R(\psi)$, which has coefficient of order $O(n^{-1/2})$, the $z^4$ and $z^6$ terms in %the marginal posterior density of $Z(\psi)$, which have coefficients of order $O(n^{-1})$, would yield a term of %order $O(n^{-3/2})$. When calculating the contribution to the expectation from the $\{Z(\psi)\}^3$ term in %$R(\psi)$, which has coefficient of order $O(n^{-1})$, the $z$ and $z^3$ terms in the marginal posterior density %of $Z(\psi)$, which have coefficients of order $O(n^{-1/2})$, would yield terms of order $O(n^{-3/2})$, while the %$z^4$ and $z^6$ terms in the marginal posterior density of $Z(\psi)$ would yield 0. Similar calculations, to %error of order $O(n^{-3/2})$, were given by DiCiccio \& Stern (1993).

\subsection*{Higher-order analysis of the Bayesian mean}

Since the focus here is the frequentist mean of $\mu_B$, it is appropriate to describe error terms in expansions by using the $O_p(\cdot)$ notation instead of the $O(\cdot)$ notation.
To simplify subsequent calculations, $\mu_B$ is conveniently written in the form $\mu_B=f(\hat\theta)+O_p(n^{-3/2})$, where
\begin{equation}
f(\theta)=-{\textstyle{1\over 2}} H L_{rst} L^{r1} L^{rs}
-{\textstyle{1\over 6}} H^3 L_{rst} L^{r1} L^{s1} L^{t1}
+ H{{\partial\log\pi(\theta)}\over{\partial\theta^r}} L^{r1}
\end{equation}
and $H=(-L^{11})^{-1/2}$. Recall that $f(\theta)$ is $O_p(n^{-1/2})$. If the prior density $\pi(\theta)$ satisfies the Welch--Peers condition then $\mu_B=\mu_F+O_p(n^{-1})$, and hence, ${ E}(\mu_B)=\mu_F+O(n^{-1})$. We establish here that ${E}(\mu_B)=\mu_F+O(n^{-3/2})$.

Note that, by Taylor expansion about $\theta$,
$$
f(\hat\theta)=f(\theta)+f_{r}(\theta)(\hat\theta^r-\theta^r)+O_p(n^{-3/2})
=f(\theta)-f_{r}(\theta)\lambda^{rs}l_s+O_p(n^{-3/2}),
$$
since $\hat\theta^r-\theta^r=\lambda^{rs}l_s+O_p(n^{-1})$. Now, $f_{r}(\theta)$ is $O_p(n^{-1/2})$, and it is a function of the $L$'s. Let $f^{\lambda}_r(\theta)$ be the quantity obtained when each of the $L$'s in $f_r(\theta)$ is replaced by its corresponding $\lambda$, so $f^{\lambda}_r(\theta)$ is a nonrandom quantity depending on $\theta$. Then
$f_r(\theta)=f^{\lambda}_r(\theta)+O_p(n^{-1}),$ and hence
${ E}\{f_{r}(\theta)\lambda^{rs}l_s\}={ E}\{f^{\lambda}_r(\theta)\lambda^{rs}l_s+O_p(n^{-3/2})\}=O(n^{-3/2}).$
It follows that
${ E}(\mu_B)={E}\{f(\theta)\}+O(n^{-3/2}),$
so it is required to show that ${E}\{f(\theta)\}=\mu_F+O(n^{-3/2})$.

We have
$L^{rs}=\lambda^{rs}-\lambda^{rt}\lambda^{su}l_{tu}+O_p(n^{-2})$,
so that
$$H=(-L^{11})^{-1/2}=\eta-{\textstyle{1\over 2}}\eta^3\lambda^{r1}\lambda^{s1}l_{rs}+O_p(n^{-1/2}),\quad
H^3=\eta^3-{\textstyle{3\over 2}}\eta^3\lambda^{r1}\lambda^{s1}l_{rs}+O_p(n^{-1/2}).$$

Consider first the final term of $f(\theta)$; the other terms can be handled similarly. It follows from the preceding equations that
\begin{eqnarray*}
H{{\partial\log\pi(\theta)}\over{\partial\theta^r}} L^{r1}&=&
\eta{{\partial\log\pi(\theta)}\over{\partial\theta^r}}\lambda^{r1}
-{\textstyle{1\over 2}}\eta^3{{\partial\log\pi(\theta)}\over{\partial\theta^r}}
\lambda^{r1}\lambda^{s1}\lambda^{t1}l_{st}\\
&&\qquad\qquad-\eta{{\partial\log\pi(\theta)}\over{\partial\theta^r}}\lambda^{rs}\lambda^{t1}l_{st}+O_p(n^{-3/2}).
\end{eqnarray*}
Hence,
\begin{eqnarray*}
{ E}\biggl\{H{{\partial\log\pi(\theta)}\over{\partial\theta^r}} L^{r1}\biggr\}&=&
\eta{{\partial\log\pi(\theta)}\over{\partial\theta^r}}\lambda^{r1}+O(n^{-3/2})\\
&=&\eta\lambda_{rs/t}\lambda^{r1}\lambda^{st}
+\textstyle{1\over 2}\eta^3\lambda_{rs/t}\lambda^{r1}\lambda^{s1}\lambda^{t1}+O(n^{-3/2}),
\end{eqnarray*}
by virtue of the Welch--Peers condition (1). The other terms in $f(\theta)$ have
\begin{eqnarray*}
{ E}(-{\textstyle{1\over 2}} H L_{rst} L^{r1} L^{rs})&=&
-{\textstyle{1\over 2}} \eta \lambda_{rst} \lambda^{r1} \lambda^{rs}+O(n^{-3/2}),\\
{ E}(-{\textstyle{1\over 6}} H^3 L_{rst} L^{r1} L^{s1} L^{t1})&=&
-{\textstyle{1\over 6}} \eta^3 \lambda_{rst} \lambda^{r1} \lambda^{s1} \lambda^{t1}+O(n^{-3/2}),
\end{eqnarray*}
and combining these expressions yields the desired result, namely
${ E}\{f(\theta)\}=\mu_F+O(n^{-3/2}).$

\section*{References}

\noindent
\begin{enumerate}
\item Barndorff-Nielsen, O.\ E.\ (1986).\ Inference on full or partial parameters based on the standardized signed log likelihood ratio.\ {\it Biometrika} {\bf 73}, 307-22.
\item Barndorff-Nielsen, O.\ E.\ \& Cox, D.\ R.\  (1994).\ {\it Inference and Asymptotics}.\ London: Chapman \& Hall.
\item
Datta, G.\ S.\ \& Mukerjee, R.\ (2004). \emph{Probability Matching Priors: Higher Order Asymptotics}. New York: Springer.
\item
DiCiccio, T.\ J.\ \& Martin, M.\ A.\  (1993).\ Simple modifications for signed roots of likelihood
ratio statistics.\  {\it J.R. Statist. Soc}.\ B {\bf 55}, 305-16.
\item
DiCiccio, T.\ J.\ \& Stern, S.\ E.\  (1993).\ On Bartlett adjustments for approximate Bayesian inference.\ {\it Biometrika} {\bf 80}, 731-40.
\item 
DiCiccio, T.\ J.\ \& Stern, S.\ E.\  (1994a).\ Constructing approximately standard normal pivots from signed roots of adjusted likelihood ratio statistics.\
{\it Scand. J. Statist.}, {\bf 21}, 447-60.
\item 
DiCiccio, T.\ J.\ \& Stern, S.\ E.\  (1994b).\ Frequentist and Bayesian Bartlett correction of test statistics based on adjusted profile likelihoods.\ {\it J.R. Statist. Soc}.\ B {\bf 56}, 397-408.
\item 
DiCiccio, T.\ J.\ \& Young, G.\ A.\  (2010). Objective Bayes and conditional inference in exponential families.\ {\it Biometrika} {\bf 97}, 497-504.
\item
Hall, P.\  (1992).\ \emph{The Bootstrap and Edgeworth Expansion}. New York: Springer.
\item
Lawless, J.\ F.\ (1982). \emph{Statistical Models and Methods for Lifetime Data}. New York: Wiley.
\item
Mukerjee, R.\ \& Dey, D.\ K.\  (1993).\ Frequentist validity of posterior quantiles in the presence of a nuisance parameter: Higher order asymptotics.{\it Biometrika}, {\bf 80}, 499-505.
\item 
Mukerjee, R.\ \& Ghosh, M.\  (1997).\ Second-order probability matching priors. {\it Biometrika}, {\bf 84}, 970-5.
\item
Nicolau, A.\ (1993). Bayesian intervals with good frequentist behavior in the presence of a nuisance parameter: higher order asymptotics. {\it J. R. Statist. Soc.} B {\bf 55}, 377-90.
\item
Peers, H.\ W.\  (1965).\ On confidence sets and Bayesian probability points in the case of several parameters.\ {\it J. Roy. Statist. Soc. B}, {\bf 27}, 9-16.
\item
Tibshirani, R.\  (1989).\ Noninformative priors for one parameter of many.\ {\it Biometrika} {\bf 76}, 604-8.
\item
Tierney, L.\ \& Kadane, J.\ B.\  (1986).\ Accurate approximations for posterior moments and marginal densities.\ {\it J. Amer. Statist. Assoc.} {\bf 81}, 82-6.
\item
Welch, B.\ L.\  \& Peers, H.\ W.\  (1963). On formulae for confidence points based on integrals of weighted likelihoods.\ {\it J. Roy. Statist. Soc B}, {\bf 62}, 159-80.
\end{enumerate}

\end{document}